\documentclass{amsart}

\usepackage{amssymb}
\usepackage{url}

\usepackage{graphicx}

\newtheorem{theorem}{Theorem}[section]
\newtheorem{lemma}[theorem]{Lemma}
\newtheorem{corollary}[theorem]{Corollary}
\newtheorem{conjecture}[theorem]{Conjecture}

\newtheorem{axiom}[theorem]{Axiom}

\newtheorem{example}[theorem]{Example}
\newtheorem{definition}[theorem]{Definition}

\def\st {\ast}
\def\difst {\divideontimes}
\def\capst {\widehat{\cap}}
\def\lch {\triangleleft}
\def\rch {\triangleright}
\def\inch {\in}

\def\totdif {\not\equiv}

\begin{document}

%%%%%%%%%%%%%%%%%%%%%%%%%
% Subject classification
%%%%%%%%%%%%%%%%%%%%%%%%%

% Provide an AMS subject classification with one or two primary classification
% numbers and, optionally, one or more secondary classification numbers.
% Use the following format:  "Primary 42B25. Secondary 42B60, 20E26"

\subjclass{46SXX}

%%%%%%%%%
% Title
%%%%%%%%%

% Title, in lower case, with no explicit linebreaks (\\).  If the title
% is too long to be used as a running head, add a short version of the
% title in brackets, as in \title[shorttitle]{fulltitle}.

\title{$\sigma$-Relations, $\sigma$-functions and $\sigma$-antifunctions.}

%%%%%%%%%%%%%%%%%%%%%%%%%%%%%%
% Author names and addresses
%%%%%%%%%%%%%%%%%%%%%%%%%%%%%%

% Provide one separate \author{...} \address{...} \email{....} entry for each
% author, i.e., do not combine multiple authors.  Separate address lines by double
% slashes.  Do not attach footnotes to author  names. (For acknowledgements use
% the "\thanks" construct below.)
%

\author{Ivan Gatica Araus}

\address{ Becario MAE-AECI}

\address{ Department of Mathematical Analysis, University of Sevilla, St.
Tarfia s/n, Sevilla, SPAIN}

\address{ Department of Mathematics, University Andr$\acute{e}$s Bello, Los Fresnos 52, Vi$\tilde{n}$a del mar, CHILE}

\email{igatica@us.es}

%%%%%%%%%%%%%%%%%%%%
% Acknowledgements
%%%%%%%%%%%%%%%%%%%

% Use \thanks for acknowledgements as footnotes to the title page.
% (Note that footnotes inside \author or \title macros are not
% allowed.)
%
% In case of multiple author papers, phrase the acknowledgement to
% say "The first author was supported by ...  The second author was
% supported by ..."

\thanks{}

%%%%%%%%%%%%%
% Abstract
%%%%%%%%%%%%%
%
% Abstracts should not contain macros (so that they can be processed independently
% of the paper.) Avoid displayed math and references in the abstract.

\begin{abstract}
In this article we develop the concepts of $\sigma$-relation and $\sigma$-function, following the same steps as in Set Theory. First we define the concept of ordered pair and then we build the Cartesian Product of $\sigma$-sets so that we can define the concepts of $\sigma$-relation and $\sigma$-function.

Now, as in $\sigma$-Set Theory there exist the concepts of $\sigma$-antielement and $\sigma$-antiset, we can build the new concepts of $\sigma$-antifunction, antidentity and antinverse. Finally, in the case that a $\sigma$-function $f:A\rightarrow B$ is bijective and there exist $A^{\st}$ and $B^{\st}$ $\sigma$-antiset of $A$ and $B$, we get 16 different $\sigma$-functions which are related in a diagram of $\sigma$-functions.
\end{abstract}

\maketitle

%%%%%%%%%%%%%%%%%%%%%%%%%%%%%%%%%%%%%%%%%%%%%%%%%%%%%%%%%%%%%%%%%%%%%%%%%
% end Topmatter
%%%%%%%%%%%%%%%%%%%%%%%%%%%%%%%%%%%%%%%%%%%%%%%%%%%%%%%%%%%%%%%%%%%%%%%%%

%%%%%%%%%%%%%%%%%%%%%%%%%%%%%%%%%%%%%%%%%%%%%%%%%%%%%%%%%%%%%%%%%%%%%%%%%
% body of paper
%%%%%%%%%%%%%%%%%%%%%%%%%%%%%%%%%%%%%%%%%%%%%%%%%%%%%%%%%%%%%%%%%%%%%%%%%

\begin{section}{Introduction}
The context in which we developed the definitions and results will be the $\sigma$-Set Theory (see \cite{Gatica}). In this sense we consider the following axiom system.

\begin{axiom}{\textbf{(Empty $\sigma$-set).}}\label{axiom empty set}
There exists a $\sigma$-set which has no $\sigma$-elements, that is
$$(\exists X)(\forall x)(x\notin X).$$
\end{axiom}

\begin{axiom}{\textbf{(Extensionality).}}\label{axiom of extensionality}
For all $\sigma$-classes $\hat{X}$ and $\hat{Y}$, if $\hat{X}$ and $\hat{Y}$ have the same
$\sigma$-elements, then $\hat{X}$ and $\hat{Y}$ are equal, that is
$$(\forall \hat{X},\hat{Y})[(\forall z)(z\in \hat{X} \leftrightarrow z\in \hat{Y})\rightarrow \hat{X}=\hat{Y}].$$
\end{axiom}

\begin{axiom}{\textbf{(Creation of $\sigma$-Class).}}\label{axiom of creation of class}
We consider an atomic formula $\Phi(x)$ (where $\hat{Y}$ is not free). Then there exists the classes $\hat{Y}$ of all $\sigma$-sets that satisfies $\Phi(x)$, that is
$$(\exists\hat{Y})(x\in \hat{Y} \leftrightarrow \Phi(x)),$$
with $\Phi(x)$ a atomic formula where $\hat{Y}$ is not free.
\end{axiom}

\begin{axiom}{\textbf{(Scheme of Replacement).}}\label{axiom of replacement}
The image of a $\sigma$-set under a normal functional formula $\Phi$ is a
$\sigma$-set.
\end{axiom}

\begin{axiom}{\textbf{(Pair).}}\label{axiom of pairs}
For all $X$ and $Y$ $\sigma$-sets there exists a $\sigma$-set $Z$, called fusion
of pairs of $X$ and $Y$, that satisfy one and only one of the
following conditions:
\begin{description}
  \item[(a)] $Z$ contains exactly $X$ and $Y$,
  \item[(b)] $Z$ is equal to the empty $\sigma$-set,
\end{description}
that is
$$(\forall X,Y)(\exists Z)(\forall a)[(a\in Z \leftrightarrow a=X \vee a=Y ) \underline{\vee}(a\notin Z)].$$
\end{axiom}

\begin{axiom}{\textbf{(Weak Regularity).}}\label{axiom of w-regularity}
For all $\sigma$-set $X$, for all $\lch x,\ldots,w\rch\in CH(X)$ we have that $X\not\inch\lch x,\ldots,w\rch$, that
is
$$(\forall X)(\forall\lch x,\ldots,w\rch\in CH(X))(X\not\inch\lch x,\ldots,w\rch).$$
\end{axiom}

\begin{axiom}{\textbf{(non $\epsilon$-Bounded $\sigma$-Set).}}\label{axiom non bounded set}
There exists a non $\epsilon$-bounded $\sigma$-set, that is
$$(\exists X)(\exists y)[(y\in X)\wedge(\min(X)=\emptyset\vee\max(X)=\emptyset)].$$
\end{axiom}

\begin{axiom}{\textbf{(Weak Choice).}}\label{axiom of weak choice}
If $\hat{X}$ be a $\sigma$-class of $\sigma$-sets, then we can choose a singleton $Y$ whose unique $\sigma$-element come from $\hat{X}$, that is
$$(\forall\hat{X})(\forall x)(x\in\hat{X}\rightarrow (\exists Y)(Y=\{x\})).$$
\end{axiom}

\begin{axiom}{\textbf{($\epsilon$-Linear $\sigma$-set).}}\label{axiom e-linear set}
There exist $\sigma$-set $X$ such that $X$ has the linear $\epsilon$-root property, that is
$$(\exists X)(\exists y)(y\in X\wedge X\in LR).$$
\end{axiom}

\begin{axiom}{\textbf{(One and One$^{\st}$ $\sigma$-set).}}\label{axiom one set}
For all $\epsilon$-linear singleton, there exists a $\epsilon$-linear singleton $Y$ such that $X$ is totally different from $Y$, that is
$$(\forall X\in SG\cap LR)(\exists Y\in SG\cap LR)(X\totdif Y).$$
\end{axiom}

\begin{axiom}{\textbf{(Completeness (A))}.}\label{axiom of completeness a}
If $X$ and $Y$ are $\sigma$-sets, then
$$\{X\}\cup\{Y\}=\{X,Y\},$$
if and only if $X$ and $Y$ satisfy one of the following conditions:
\begin{description}
  \item[(a)] $\min(X,Y)\neq |1\vee 1^{\st}|\wedge \min(X,Y)\neq |1^{\st}\vee 1|.$
  \item[(b)] $\neg(X\totdif Y).$
  \item[(c)] $(\exists w\in X)[w\notin \min(X)\wedge \neg \Psi(z,w,a,Y)].$
  \item[(d)] $(\exists w\in Y)[w\notin \min(Y)\wedge \neg \Psi(z,w,a,X)].$
\end{description}
\end{axiom}

\begin{axiom}{\textbf{(Completeness (B)).}}\label{axiom of completeness b}
If $X$ and $Y$ are $\sigma$-sets, then
$$\{X\}\cup\{Y\}=\emptyset,$$
if and only if $X$ and $Y$ satisfy the
following conditions:

\begin{description}
  \item[(a)] $\min(X,Y)=|1\wedge 1^{\st}|\vee \min(X,Y)=|1^{\st}\wedge 1|$;
  \item[(b)] $X\totdif Y;$
  \item[(c)] $(\forall z)(z\in X\wedge z\notin \min(X))\rightarrow \Psi(z,w,a,Y))$;
  \item[(d)] $(\forall z)(z\in Y\wedge z\notin \min(Y))\rightarrow \Psi(z,w,a,X))$.
\end{description}
\end{axiom}

\begin{axiom}{\textbf{(Exclusion).}}\label{axiom of exclusion}
For all $\sigma$-sets $X,Y,Z$, if $Y$ and $Z$ are $\sigma$-elements of $X$ then the fusion of pairs of $Y$ and $Z$ contains exactly $Y$ and $Z$, that is
$$(\forall X,Y,Z)(Y,Z\in X\rightarrow \{Y\}\cup\{Z\}=\{Y,Z\}).$$
\end{axiom}

\begin{axiom}{\textbf{(Power $\sigma$-set).}}\label{axiom of power set}
For all $\sigma$-set $X$  there exists a $\sigma$-set $Y$, called power of $X$, whose $\sigma$-elements are exactly the $\sigma$-subsets of $X$, that is
$$(\forall X)(\exists Y)(\forall z)(z\in Y\leftrightarrow z\subseteq X).$$
\end{axiom}

\begin{axiom}{\textbf{(Fusion).}}\label{axiom of fusion}
For all $\sigma$-sets $X$ and $Y$, there exists a $\sigma$-set $Z$, called fusion of all $\sigma$-elements of $X$ and $Y$, such that $Z$ contains $\sigma$-elements of the $\sigma$-elements of $X$ or $Y$, that is
$$(\forall X,Y)(\exists Z)(\forall b)(b \in Z \rightarrow (\exists z)[(z \in X\vee z\in Y)\wedge(b\in z)]).$$
\end{axiom}

\begin{axiom}{\textbf{(Generated $\sigma$-set).}}\label{axiom of generated set}
For all $\sigma$-sets $X$ and $Y$ there exists a $\sigma$-set, called the $\sigma$-set generated by $X$ and $Y$, whose $\sigma$-elements are exactly the fusion of the $\sigma$-subsets of $X$ with the $\sigma$-subsets of $Y$, that is
$$(\forall X,Y)(\exists Z)(\forall a)(a \in Z \leftrightarrow(\exists A\in 2^{X})(\exists B\in 2^{Y})(a=A\cup B)).$$
\end{axiom}

Where:
\begin{itemize}
  \item $\lch x,\ldots,z\rch\in CH(X):= x\in \cdots\in z\in X.$
  \item $X\notin \lch x,\ldots,z\rch:= X\neq x\wedge \cdots\wedge X\neq z.$
  \item $LR=\{X: X \textrm{ has the linear $\epsilon$-root property }\}.$
  \item $SG=\{X: X \textrm{ is a singleton }\}.$
	\item $\Psi(z,w,a,x):=(\exists ! w)(\{z\}\cup\{w\}=\emptyset) \wedge (\forall a)(\{z\}\cup\{a\}=\emptyset\rightarrow a\in x).$
	\item $1$ is the One $\sigma$-set and $1^{\st}$ is the One$^{\st}$ $\sigma$-set.
	\item $min(X)=\{y\in X: y \textrm{ is an $\epsilon$-minimal $\sigma$-element of } X\}$.
	\item $\min(X,Y)\neq |1\vee 1^{\st}|:= \min(X)\neq 1\vee \min(Y)\neq 1^{\st}.$
	\item $\min(X,Y)=|1\wedge 1^{\st}|:= \min(X)=1\wedge \min(Y)= 1^{\st}.$	
\end{itemize}

\begin{figure}
\centering
\includegraphics[width=1\textwidth]{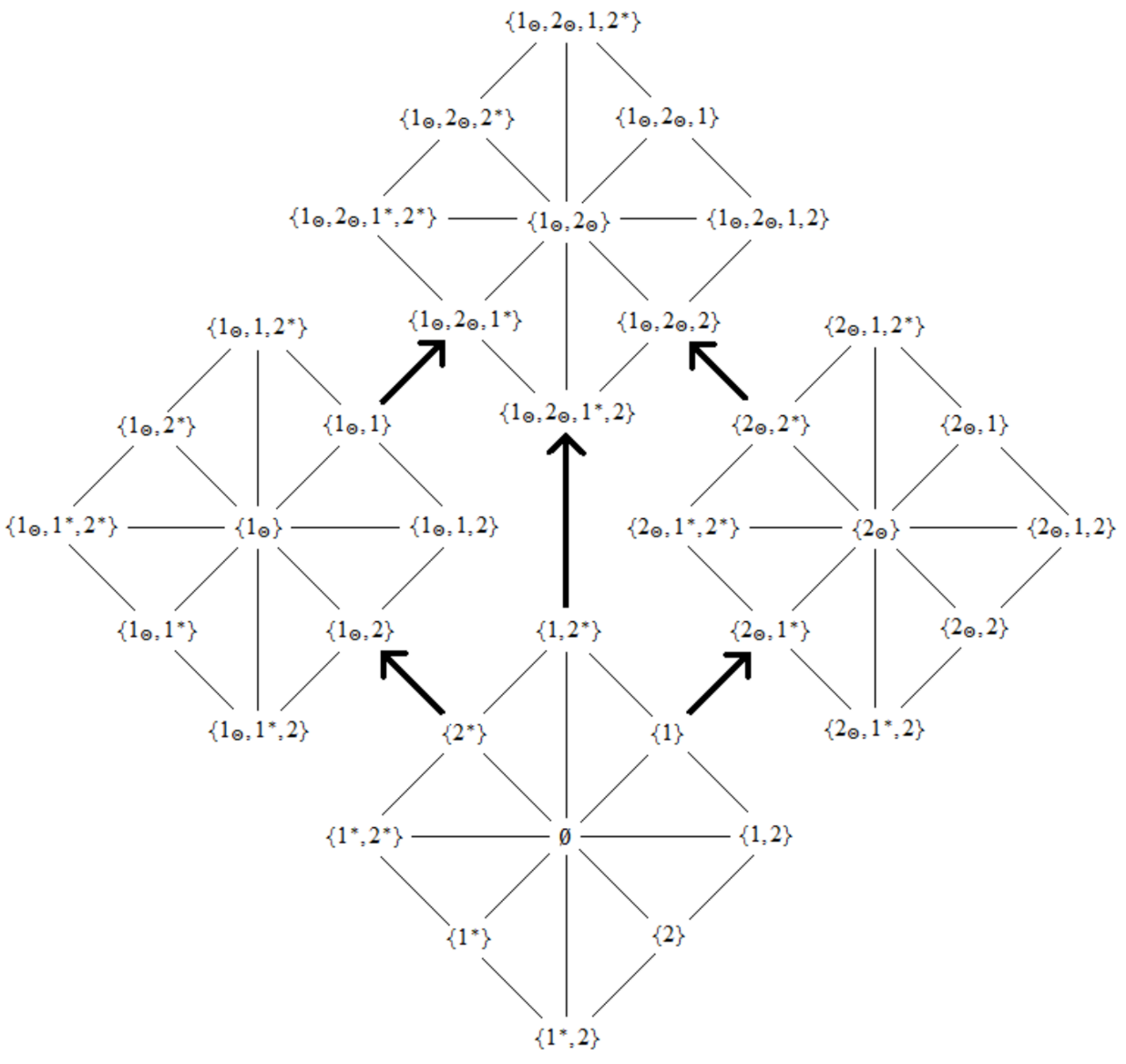}
\caption{Generated Space by $\{1_{\Theta},2_{\Theta},1,2\}$ and $\{1_{\Theta},2_{\Theta},1^{\st},2^{\st}\}$}
\label{generated space}
\end{figure}

Also, we remember that the generic object of $\sigma$-Set Theory is called $\sigma$-class. However, in what follows we consider only the properties of $\sigma$-sets. In this sense when we write $x\in A$, it should be understood that $x$ and $A$ are $\sigma$-sets, where $x$ is a $\sigma$-element of $A$.

Now, we present the following definitions and results introduced by I. Gatica in \cite{Gatica} which are necessary for the development of this paper.

\begin{definition}\label{def 1 antielement and antiset}
Let $A$ and $B$ be $\sigma$-sets. Then we have that
\begin{enumerate}
	\item If $\{A\}\cup\{B\}=\emptyset$, then $B$ \textbf{is called the $\sigma$-antielement of } $A$.
	\item If $A\cup B=\emptyset$, then $B$ \textbf{is called the $\sigma$-antiset of } $A$.
\end{enumerate}
\end{definition}

To denote the $\sigma$-antielements and $\sigma$-antisets we use the following notation:
\begin{itemize}
	\item Let $x\in A$, then we use $x^{\st}$ to denote the $\sigma$-antielement of $x$.
	\item Let $A$ be a $\sigma$-set, then we use $A^{\st}$ to denote the $\sigma$-antiset of $A$.
\end{itemize}

We observe that Gatica in \cite{Gatica} uses $A^{\star}$ in order to denote $\sigma$-antisets. However, in order to unify notation, we use $A^{\st}$.

\begin{definition}\label{def 1 int and diff antiset}
Let $A$ and $B$ be $\sigma$-sets. We define:
\begin{enumerate}
  \item $A\cap B=\{x\in A: x\in B\};$
  \item $A-B=\{x\in A: x\notin B\};$
  \item $A\capst B :=\{x\in A: x^{\st}\in B\};$
  \item $A\difst B:=A-(A\capst B);$
  \item $A\cup B=\{x: (x\in A\difst B)\vee(x\in B\difst A)\}.$
\end{enumerate}
\end{definition}

We observe that if $A\capst B=\emptyset$ and $B\capst A=\emptyset$ then 
$$A\cup B=\{x: x\in A\vee x\in B\}.$$
Thus, in this case the fusion coincides with the definition of union in a standard set theory. 

\begin{example}\label{example fusion equiv union}
Let $X$ be a nonempty $\sigma$-set and $2^{X}$ the power $\sigma$-set of $X$. Then for all $A,B\in 2^{X}$ we have that  
$$A\cup B=\{x: x\in A\vee x\in B\},$$
because $2^{X}$ is a $\sigma$-antielement free $\sigma$-set (see \cite{Gatica}, Definition 3.44 and Theorem 3.48).   
\end{example}

\begin{definition}\label{def 1 generated space by A and B}
Let $A$ and $B$ be $\sigma$-sets. The \textbf{Generated space by $A$
and $B$} is given by
$$\langle 2^{A},2^{B}\rangle=\{x\cup y: x\in 2^{A} \wedge y\in 2^{B}\}.$$
\end{definition}

We observe that in general $2^{A\cup B}\neq \langle 2^{A},2^{B}\rangle.$ Consider $X=\{1_{\Theta}2^{\st}\}$ and $Y=\{1_{\Theta},2\}$, then $X\cup Y=\{1_{\Theta}\}$. Therefore $2^{X\cup Y}=\{\emptyset, \{1_{\Theta}\}\}$ and $\langle 2^{X},2^{Y}\rangle=\{\emptyset, \{1_{\Theta}\}, \{2\}, \{2^{\st}\}, \{1_{\Theta},2\},\{1_{\Theta},2^{\st}\} \}.$

On the other hand, if we consider $A=\{1_{\Theta},2_{\Theta}\}$ and $B=\{1,2\}$ we obtain that the generated space by $A\cup B$ 
and $A\cup B^{\st}$ is the following:

\

$\langle 2^{A\cup B},2^{A\cup B^{\st}}\rangle=$ $\{\emptyset, \{1\},$ $\{1^{\st}\},$ $\{2\},$ $\{2^{\st}\},$ $\{1_{\Theta}\},$ $\{2_{\Theta}\},$ $\{1^{\st},2\},$ $\{1,2^{\st}\},$ $\{1,2\},$ $\{1^{\st},2^{\st}\},$ $\{1_{\Theta},1\},$ $\{1_{\Theta},1^{\st}\},$ $\{1_{\Theta},2\},$ $\{1_{\Theta},2^{\st}\},$ $\{2_{\Theta},1\},$ $\{2_{\Theta},1^{\st}\},$ $\{2_{\Theta},2\},$ $\{2_{\Theta},2^{\st}\},$ $\{1_{\Theta},2_{\Theta}\},$
$\{1_{\Theta},1^{\st},2\},$ $\{1_{\Theta},1,2^{\st}\},$ $\{1_{\Theta},1,2\},$ $\{1_{\Theta},1^{\st},2^{\st}\},$ $\{2_{\Theta},1^{\st},2\},$ $\{2_{\Theta},1,2^{\st}\},$ $\{2_{\Theta},1,2\},$ \ $\{2_{\Theta},1^{\st},2^{\st}\},$ \ $\{1_{\Theta},2_{\Theta},1\},$ $\{1_{\Theta},2_{\Theta},1^{\st}\},$ \ $\{1_{\Theta},2_{\Theta},2\},$ $\{1_{\Theta},2_{\Theta},2^{\st}\},$ $\{1_{\Theta},2_{\Theta},1^{\st},2\},$ $\{1_{\Theta},2_{\Theta},1,2^{\st}\},$ $\{1_{\Theta},2_{\Theta},1,2\},$ $\{1_{\Theta},2_{\Theta},1^{\st},2^{\st}\}\}.$

\

See Figure \ref{generated space}. Also if $A=\{1,2,3\}$ and $A^{*}=\{1^{\st},2^{\st},3^{\st}\}$, then

\

$\langle 2^{A},2^{A^{\st}}\rangle=3^{A}$ $=\{\emptyset,$ $\{1\},$ $\{2\},$ $\{3\},$ $\{1^{\st}\},$ $\{2^{\st}\},$ $\{3^{\st}\},$ $\{1,2\},$ $\{1,3\},$ $\{2,3\},$ $\{1^{\st},2\},$ $\{1^{\st},3\},$ $\{2^{\st},3\},$ $\{1^{\st},2^{\st}\},$ $\{1^{\st},3^{\st}\},$ $\{2^{\st},3^{\st}\},$ $\{1,2^{\st}\},$ $\{1,3^{\st}\},$ $\{2,3^{\st}\},$ $\{1,2,3\},$ \ $\{1^{\st},2,3\},$ \ $\{1,2^{\st},3\},$ \ $\{1,2,3^{\st}\},$ $\{1^{\st},2^{\st},3\},$ $\{1^{\st},2,3^{\st}\},$ $\{1,2^{\st},3^{\st}\},$ $\{1^{\st},2^{\st},3^{\st}\} \}.$

\

Now, in Figure \ref{integer space}, we present the patterns of containments of the $3^{X}$. We observe that the $\sigma$-elements $\{1^{\st},2,3^{\st}\}$ and $\{1,2^{\st},3\}$ can be represented in three dimensions as one of
the vertexes of the pyramids
$$\triangle_{1}:=\{\{1,3\},\{2^{\st},3\},\{1,2^{\st}\}, \{1,2^{\st},3\}\}$$
and
$$\triangle_{2}:= \{\{1^{\st},3^{\st}\},\{2,3^{\st}\},\{1^{\st},2\},\{1^{\st},2,3^{\st}\}\}.$$

\begin{figure}
\centering
\includegraphics[width=0.75\textwidth]{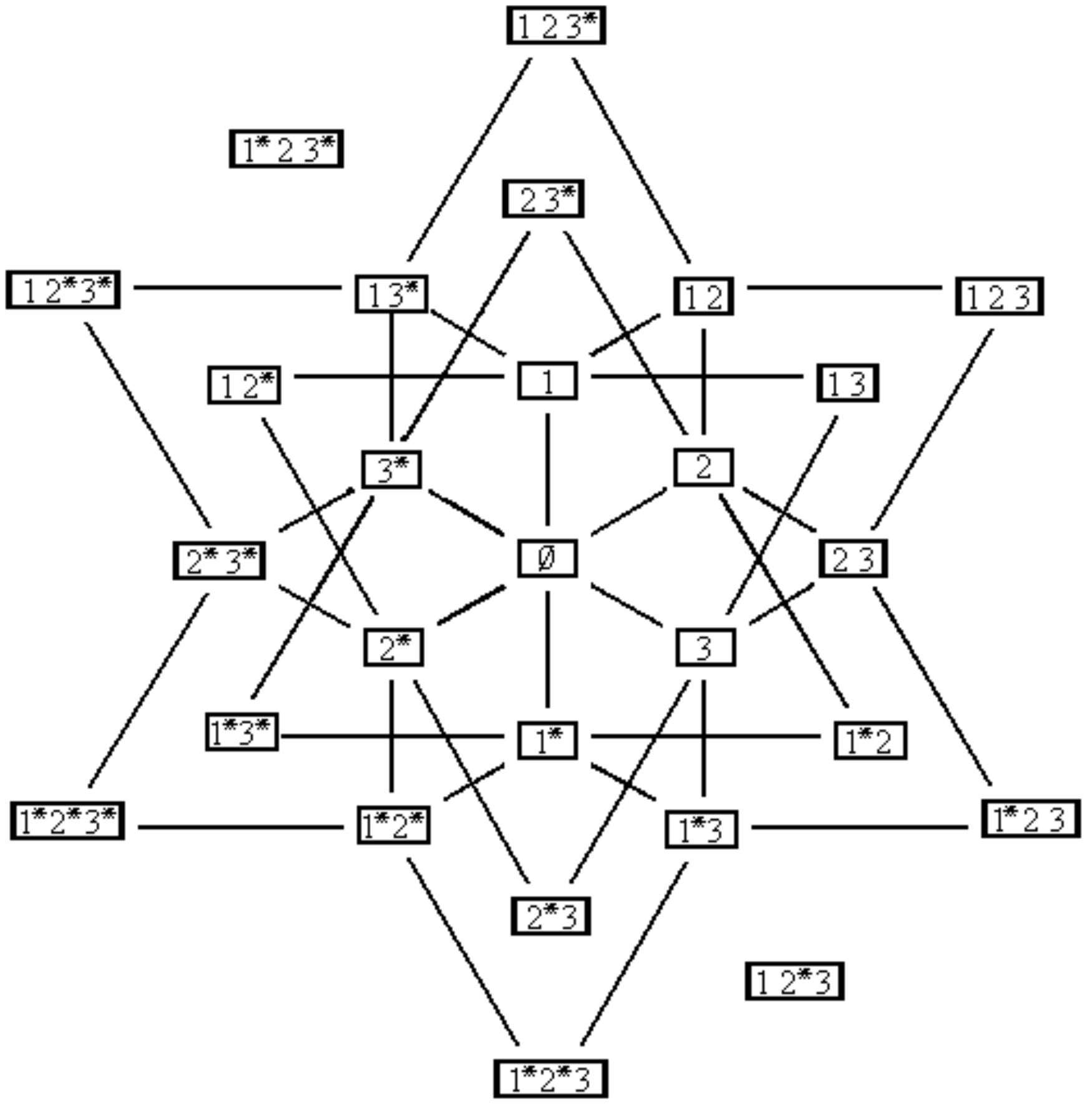}
\caption{Integer Space $3^{\{1,2,3\}}$}
\label{integer space}
\end{figure}

\newpage

\begin{theorem}\label{theo 1 fusion of pairs}
(see [\cite{Gatica}, Theorem 3.29]). If $A$ is a $\sigma$-set, then
\begin{description}
  \item[(a)] $\{\emptyset\}\cup\{A\}=\{\emptyset,A\}.$
  \item[(b)] $\{\alpha\}\cup\{A\}=\{\alpha,A\}.$
  \item[(c)] $\{\beta\}\cup\{A\}=\{\beta,A\}.$
\end{description}
\end{theorem}

\begin{theorem}\label{theo 1 oper pairs is conmutative}
(see [\cite{Gatica}, Theorem 3.32]). Let $A$ and $B$ be $\sigma$-set. Then $\{A\}\cup\{B\}=\emptyset$ if and only if $\{B\}\cup\{A\}=\emptyset$.
\end{theorem}

\begin{theorem}\label{theo 1 antielement is unique}
(see [\cite{Gatica}, Theorem 3.34]). Let $A$ be a $\sigma$-set. If there exists the $\sigma$-antielement of $A$, then it is unique.
\end{theorem}

\begin{theorem}\label{theo 1 antiset is unique}
(see [\cite{Gatica}, Theorem 3.58]). Let $A$ be a $\sigma$-set. If there exists $A^{\st}$, the $\sigma$-antiset of $A$, then $A^{\st}$ is unique.
\end{theorem}

\end{section}

\newpage

\begin{section}{Order Pairs}

Following the same steps as in the Set Theory, we define the ordered pairs. Next we define the Cartesian product of the $\sigma$-sets. Therefore we can define the concept of $\sigma$-relations, $\sigma$-functions and $\sigma$-antifunctions.

\begin{definition}\label{def 1 order pairs}
Let $X$ and $Y$ be $\sigma$-sets. The ordered pair of $X$ and $Y$ is defined by
$$(X,Y):=\{\{X\}\}\cup\{\{X\}\cup\{Y\}\}.$$
\end{definition}

We observe that for all $X$ and $Y$ $\sigma$-sets, the ordered pair $(X,Y)$ is a
$\sigma$-set.

\begin{lemma}\label{lemma 1 order pairs}
Let $X$ and $Y$ be $\sigma$-set. Then $(X,Y)=\{\{X\},\{X\}\cup\{Y\}\}.$
\end{lemma}

\begin{proof}
Consider $A=\{X\}$ and $B=\{X\}\cup\{Y\}$. By Axiom 4 (Pairs) we
have that $B=\{X,Y\}$ or $B=\emptyset$. If $B=\{X,Y\}$ then $B$ is
not totally different from $A$. Therefore by the condition (b) of
Axiom 7 (Completeness (A)) we obtain that $(X,Y)=\{\{X\},\{X,Y\}\}$.
Now, if $B=\emptyset$ then by Theorem \ref{theo 1 fusion of pairs},
we obtain that $(X,Y)=\{\{X\},\emptyset\}$. Therefore
$(X,Y)=\{\{X\},\{X\}\cup\{Y\}\}.$
\end{proof}

\begin{theorem}\label{theo 1 order pairs}
Let $X$, $Y$, $Z$ and $W$ be $\sigma$-sets. Then $(X,Y)=(Z,W)$ if and only if $(X=Z\wedge Y=W).$
\end{theorem}

\begin{proof}

$(\rightarrow)$ Consider $(X,Y)=(Z,W)$ then we will prove that
$(X=Z\wedge Y=W)$.

\begin{description}
  \item[(case a)] Suppose that $X=\emptyset$. By Theorem \ref{theo 1 fusion of pairs}
and Lemma \ref{lemma 1 order pairs} we obtain
$(X,Y)=\{\{\emptyset\},\{\emptyset, Y\}\}$. Since $(X,Y)=(Z,W)$ then
$$\{\{\emptyset\},\{\emptyset, Y\}\}=\{\{Z\},\{Z\}\cup\{W\}\}.$$
It is clear by Axiom 4 (Pairs) that $\{Z\}\cup\{W\}=\{Z,W\}$ or
$\{Z\}\cup\{W\}=\emptyset$. Since $(X,Y)= \{\{Z\},\{Z\}\cup\{W\}\}$
then by Axiom 2 (Extensionality) we obtain that
$\{Z\}\cup\{W\}\neq\emptyset$, in consequence
$\{Z\}\cup\{W\}=\{Z,W\}$. Therefore $\{\emptyset\}=\{Z\}$ and
$\{\emptyset, Y\}=\{Z,W\}$. Finally, $Z=\emptyset$ and so $Y=W$.

  \item[(case b)] Suppose that $Y=\emptyset$. By Theorem \ref{theo 1 fusion of pairs}
and Lemma \ref{lemma 1 order pairs} we obtain
$(X,Y)=\{\{X\},\{X,\emptyset\}\}$. Now, if we use the same argument
as in (case a) we have that $X=Z$ and $W=\emptyset$.

  \item[(case c)] Suppose that $X\neq\emptyset$ and
  $Y\neq\emptyset$. By Lemma \ref{lemma 1 order pairs} we obtain
$(X,Y)=\{\{X\},\{X\}\cup\{Y\}\}$ and
$(Z,W)=\{\{Z\},\{Z\}\cup\{W\}\}$. It is clear by Axiom 4 (Pairs)
that $\{X\}\cup\{Y\}=\{X,Y\}$ or $\{X\}\cup\{Y\}=\emptyset$.

\begin{description}
  \item[(c.1)] Consider $\{X\}\cup\{Y\}=\{X,Y\}$. Now, if we use the same argument
as in (case a) we have that $X=Z$ and $Y=W$.

  \item[(c.2)] If $\{X\}\cup\{Y\}=\emptyset$ then $Y$ is the $\sigma$ antielement of $X$. Since $(X,Y)=(Z,W)$, by Axiom 2 (Extensionality) we obtain that $\{X\}=\{Z\}$ and
  $\{Z\}\cup\{W\}=\emptyset$. Finally, we obtain that $X=Z$, in consequence by Theorem \ref{theo 1 antielement is unique} $Y=W$.
\end{description}
\end{description}

$(\leftarrow)$ Consider $(X=Z\wedge Y=W)$ then by Axiom 2 (Extensional) it is clear that $(X,Y)=(Z,W)$.
\end{proof}

\begin{lemma}
Let $A$ and $B$ be $\sigma$-sets. If $x\in A$ and $y\in B$, then
$(x,y)\in 2^{\langle 2^{A},2^{B}\rangle}$.
\end{lemma}

\begin{proof}
Let $x\in A$ and $y\in B$. By Definition \ref{def 1 generated space
by A and B} we have that
$$\langle 2^{A},2^{B}\rangle=\{a\cup b: a\in 2^{A} \wedge b\in 2^{B}\}.$$
It is clear that $\{x\}\in 2^{A}$ and $\{y\}\in 2^{B}$. Therefore
$\{x\}$, $\{y\}$ and $\{x\}\cup\{y\}$ are $\sigma$-elements of
$\langle 2^{A},2^{B}\rangle$. Finally, since
$(x,y)=\{\{x\},\{x\}\cup\{y\}\}$ then $(x,y)\in 2^{\langle
2^{A},2^{B}\rangle}$.
\end{proof}

\begin{definition}
Let $A$ and $B$ be $\sigma$-sets. The \textbf{Cartesian product of
$A$ and $B$} is the $\sigma$-set of all order pairs $(x,y)$ such
that $x\in A$ and $y\in B$, that is
$$A\times B=\{(x,y): x\in A\wedge y\in B\}.$$
\end{definition}

The Cartesian product $A\times B$ is a $\sigma$-set because
$$A\times B\subset 2^{\langle 2^{A},2^{B}\rangle}.$$
\end{section}

\begin{section}{$\sigma$-Relations, $\sigma$-functions and $\sigma$-antifunctions.}

Now, we present the definition of binary $\sigma$-relations, as in
the Set Theory.

\begin{definition}
Let $A$ and $B$ be $\sigma$-sets.
\begin{description}
  \item[(a)] A \textbf{binary $\sigma$-relation} on $A\times B$ is a $\sigma$-subset $R$ of
$A\times B$.
  \item[(b)] A \textbf{binary $\sigma$-relation} on $A$ is a $\sigma$-subset $R$ of
$A\times A$.
\end{description}

If $R$ is a binary $\sigma$-relation on $A\times B$, then we also
use $R(x)=y$ for $(x,y)\in R$. As in the Set Theory we define the
\textbf{domain} of $R$ as
$$dom(R)=\{x\in A: (\exists y\in B)(R(x)=y)\},$$
and the \textbf{range} of $R$ as
$$ran(R)=\{y\in B: (\exists x\in A)(R(x)=y)\}.$$

\begin{definition}
Let $A$ and $B$ be $\sigma$-sets. A binary $\sigma$-relation $f$ on
$A\times B$ is a \textbf{$\sigma$-function} from $A$ to $B$ if for
all $x\in A$ there exists a unique $y\in B$ such that $f(x)=y$. The
unique $y$ such that $f(x)=y$ is the \textbf{value} of $f$ at $x$.
\end{definition}

Let $f$ be a $\sigma$-function on $A\times B$, then it is clear that
$dom(f)=A$ and $ran(f)\subseteq B$. Also, the $\sigma$-function from
$A$ to $B$ will be denoted by
$$f:A\rightarrow B.$$

A \textbf{binary operation on $A$} is a $\sigma$-function from
$A\times A$ to $A$ \\(i.e. $f:A\times A\rightarrow A$).

\begin{example}
Let $A$ be a $\sigma$-set. If we define $\oplus: 2^{A}\times
2^{A}\rightarrow 2^{A}$ were $\oplus(x,y)=x\cup y$, then $\oplus$ is
a \textbf{binary operation on $2^{A}$}.

In fact, let $(x,y)\in 2^{A}\times 2^{A}$ then $x,y\in 2^{A}$. Since
$x\subseteq A$ and $y\subseteq A$ then it is clear that there exist
$x\cup y\in 2^{A}$ and by Axiom 2 (Extensionality) it is unique.
\end{example}
\end{definition}

The definition of binary operation, is important for the study of
the algebraic properties of the Integer Space. In this sense, we
present the following conjecture:

\begin{conjecture}\label{conj 3 binary operation integer space}
Let $A$ be a $\sigma$-set such that there exists $A^{\st}$, the
$\sigma$-antiset of $A$. If we define $\oplus: 3^{A}\times
3^{A}\rightarrow 3^{A}$ where $\oplus(x,y)=x\cup y$, then $\oplus$ is
a \textbf{binary operation on $3^{A}$}.
\end{conjecture}

For convenience, we shall use the following notation: Let $A,B,C$ be $\sigma$-set then 
$$\delta_{A}:=\{A\},\delta_{AB}:=\{A,B\} \ \textrm{ and } \ \delta_{ABC}:=\{A,B,C\}.$$

\

For example, if we consider the integer space $3^{\{1,2,3\}}$ (see, figure \ref{integer space}) we obtain that $\oplus$ is
a binary operation on $3^{\{1,2,3\}}$.

\begin{center}
\begin{tabular}{||c||c|c|c|c|c|c|c|c||}
\hline
\hline
$\oplus$      & $\emptyset$ & $\delta_{1}$ & $\delta_{2}$    & $\delta_{3}$ & $\delta_{12}$ & $\delta_{13}$ & $\delta_{23}$ & $\delta_{123}$ \\
\hline
\hline

$\emptyset$ & $\emptyset$ & $\delta_{1}$           & $\delta_{2}$    & $\delta_{3}$ & $\delta_{12}$ & $\delta_{13}$ & $\delta_{23}$ & $\delta_{123}$ \\

\hline

$\delta_{1^{*}}$          & $\delta_{1^{*}}$          & $\emptyset$& $\delta_{1^{*}2}$ & $\delta_{1^{*}3}$ & $\delta_{2}$ & $\delta_{3}$ & $\delta_{1^{*}23}$ & $\delta_{23}$ \\

\hline

$\delta_{2^{*}}$ & $\delta_{2^{*}}$ & $\delta_{12^{*}}$ & $\emptyset$ & $\delta_{2^{*}3}$ & $\delta_{1}$ & $\delta_{12^{*}3}$ & $\delta_{3}$ & $\delta_{13}$ \\

\hline

$\delta_{3^{*}}$ & $\delta_{3^{*}}$ & $\delta_{13^{*}}$ & $\delta_{23^{*}}$ & $\emptyset$ & $\delta_{123^{*}}$ & $\delta_{1}$ & $\delta_{2}$ & $\delta_{12}$ \\

\hline

$\delta_{1^{*}2^{*}}$ & $\delta_{1^{*}2^{*}}$ & $\delta_{2^{*}}$ & $\delta_{1^{*}}$ & $\delta_{1^{*}2^{*}3}$ & $\emptyset$ & $\delta_{2^{*}3}$ & $\delta_{1^{*}3}$ & $\delta_{3}$ \\

\hline

$\delta_{1^{*}3^{*}}$ & $\delta_{1^{*}3^{*}}$ & $\delta_{3^{*}}$ & $\delta_{1^{*}23^{*}}$ & $\delta_{1^{*}}$ & $\delta_{23^{*}}$ & $\emptyset$ & $\delta_{1^{*}2}$ & $\delta_{2}$ \\

\hline

$\delta_{2^{*}3^{*}}$ & $\delta_{2^{*}3^{*}}$ & $\delta_{12^{*}3^{*}}$ & $\delta_{3^{*}}$ & $\delta_{2^{*}}$ & $\delta_{13^{*}}$ & $\delta_{12^{*}}$ & $\emptyset$ & $\delta_{1}$ \\

\hline

$\delta_{1^{*}2^{*}3^{*}}$ & $\delta_{1^{*}2^{*}3^{*}}$ & $\delta_{2^{*}3^{*}}$ & $\delta_{1^{*}3^{*}}$ & $\delta_{1^{*}2^{*}}$ & $\delta_{3^{*}}$ & $\delta_{2^{*}}$ & $\delta_{1^{*}}$ & $\emptyset$ \\
\hline
\hline
\end{tabular}
\end{center}

\

Conjecture \ref{conj 3 binary operation integer space} will be studied in future works. Now, we present the definition of
$\sigma$-antifunctions.

\begin{definition}\label{def 3 antifunction}
Let $A$, $B$ and $C$ be $\sigma$-sets. If $f: A\rightarrow B$ and
$f^{\st}: A\rightarrow C$ are $\sigma$-functions, then we say that
$f^{\st}$ is the $\sigma$-antifunction of $f$ if for all $x\in A$ we
have that $\{f(x)\}\cup\{f^{\st}(x)\}=\emptyset$ (i.e.
$(f^{\st}(x))^{\st}=f(x)$).
\end{definition}

\begin{example}\label{example 3 antifunction}
Consider the following $\sigma$-sets $A=\{1_{\Theta},2_{\Theta},3_{\Theta}\}$, $B=\{1,2,3\}$ and $C=\{1^{\st},2^{\st},4\}$.
Now, we define $f: A\rightarrow B$ such that $f(1_{\Theta})=1$, $f(2_{\Theta})=2$ and
$f(3_{\Theta})=2$ and $f^{\st}: A\rightarrow C$ such that $f(1_{\Theta})=1^{\st}$,
$f(2_{\Theta})=2^{\st}$ and $f(3_{\Theta})=2^{\st}$. It is clear that $f^{\st}$ is the
$\sigma$-antifunction of $f$. Also, we obtain that $ran(f)=\{1,2\}$
and $ran(f^{\st})=\{1^{\st},2^{\st}\}$.
\end{example}

In this sense we obtain the following Theorem.

\begin{theorem}\label{theo 3 ran(func) antset of ran(antifunc)}
Let $A$, $B$ and $C$ be $\sigma$-sets. If $f: A\rightarrow B$ is a
$\sigma$-function and $f^{\st}: A\rightarrow C$ is the
$\sigma$-antifunction of $f$, then $ran(f)$ is the $\sigma$-antiset
of $ran(f^{\st})$ (i.e. $ran(f)^{\st}=ran(f^{\st})$ and
$ran(f)=ran(f^{\st})^{\st}$ ).
\end{theorem}

\begin{proof}
Let $f: A\rightarrow B$ a $\sigma$-function and $f^{\st}:
A\rightarrow C$ the $\sigma$-antifunction of $f$. We will prove that
$ran(f)\cup ran(f^{\st})=\emptyset$. It is clear that
$$ran(f)=\{y\in B: (\exists x\in A)(f(x)=y)\},$$
$$ran(f^{\st})=\{y\in C: (\exists x\in A)(f^{\st}(x)=y)\}.$$
Now, by Definition we obtain that
$$ran(f)\cup ran(f^{\st})=\{x: (x\in ran(f)\difst ran(f^{\st}))\vee(x\in ran(f^{\st})\difst ran(f))\}.$$
Then, in order to prove that $ran(f)\cup ran(f^{\st})=\emptyset$ it is enough
to prove that
\begin{description}
  \item[(a)] $ran(f)\capst ran(f^{\st})=ran(f),$
  \item[(b)] $ran(f^{\st})\capst ran(f)=ran(f^{\st})$.
\end{description}
We will only prove (a) because the proof of (b) is similar. It is clear by Definition \ref{def 1 int and diff antiset} that $ran(f)\capst ran(f^{\st})\subseteq ran(f)$. Now, let $y\in ran(f)$ then there exists $\widehat{x}\in A$ such that $f(\widehat{x})=y$.
Since $f^{\st}$ is the $\sigma$-antifunction of $f$, then $\{f(\widehat{x})\}\cup\{f^{\st}(\widehat{x})\}=\emptyset$.
Therefore, by Theorem \ref{theo 1 antielement is unique}
$y^{\st}=f^{\st}(\widehat{x})\in ran(f^{\st})$ and so $y\in ran(f)\capst ran(f^{\st})$. Finally, $ran(f)= ran(f)\capst ran(f^{\st})$.
\end{proof}

\begin{theorem}\label{theo 3 ex anti ran(f) then  ex antifuntion}
Let $A$ and $B$ be $\sigma$-sets and $f:A\rightarrow B$ a $\sigma$-function. If there exists the $\sigma$-antiset of $ran(f)$ (i.e. $ran(f)^{\st}$), then there exists $f^{\st}$ the
$\sigma$-antifunction of $f$.
\end{theorem}

\begin{proof}
Since there exists $ran(f)^{\st}$, then we can define a $\sigma$-set $C$ such that $ran(f)^{\st}\subseteq C$. Now, we define $f^{\st}:A\rightarrow C$ such that $f^{\st}(x)=(f(x))^{\st}$. It is clear by definition that for all $x\in A$ we have that
$\{f(x)\}\cup\{f^{\st}(x)\}=\emptyset$. Therefore $f^{\st}$ is the $\sigma$-antifunction of $f$.
\end{proof}

Let $f:A\rightarrow B$ and $g:A\rightarrow C$ be $\sigma$-functions. We say that $f$ and $g$
are equal if and only if for all $x\in A$, $f(x)=g(x)$.

\begin{theorem}
Let $A$, $B$ and $C$ be $\sigma$-sets. If $f: A\rightarrow B$ is a
$\sigma$-function and $f^{\st}: A\rightarrow C$ is the
$\sigma$-antifunction of $f$, then $f^{\st}$ is unique.
\end{theorem}

\begin{proof}
This proof is obvious by Theorems \ref{theo 1 antiset is unique} and
\ref{theo 3 ran(func) antset of ran(antifunc)}.
\end{proof}

We observe that some consequences of the uniqueness of $\sigma$-antielements, $\sigma$-antisets and $\sigma$-antifunctions are the following:
\begin{itemize}
  \item If $x^{\st}$ is the $\sigma$-antielement of $x$, then $(x^{\st})^{\st}=x$.
  \item If $A^{\st}$ is the $\sigma$-antiset of $A$, then $(A^{\st})^{\st}=A$.
	\item If $f^{\st}$ is the $\sigma$-antifunction of $f$, then $(f^{\st})^{\st}=f$.
\end{itemize}

Now, following the same steps as in Set Theory we define the following: If $f:A\rightarrow B$ is a $\sigma$-function then
\begin{itemize}
  \item $f$ is a $\sigma$-function \textbf{onto} $B$ iff $B=ran(f)$.
  \item $f$ is a $\sigma$-function \textbf{one-one} iff for all
$x,y\in A$, $x\neq y \rightarrow f(x)\neq f(y).$
  \item $f$ is a $\sigma$-function \textbf{bijective} iff it is
both one-one and onto.
\end{itemize}

Also, we define

\begin{itemize}
  \item The \textbf{image} of $A$ by $f$, $f(A)=\{y\in B: (\exists x\in A)(f(x)=y)\}=ran(f)$.
  \item The \textbf{preimage} of $A$ under $f$, $f_{-1}(A)=\{x\in A: f(x)\in B\}$.
\end{itemize}

Now, it is clear that: If $f$ is bijective there is a unique
$\sigma$-function $f_{-1}:B\rightarrow A$ (called the inverse of
$f$) such that
$$(\forall x\in A)(f_{-1}(f(x))=x).$$

We have chosen the notation $f_{-1}$ to denote the inverse $\sigma$-functions for convenience.

Now, we consider $f:A\rightarrow B$ and $g:C\rightarrow D$ be
$\sigma$-functions. If $ran(g)\cap dom(f)\neq\emptyset$, then we can define the composition
of $f$ and $g$ is the $\sigma$-function $f\circ g$ with domain
$dom(f\circ g)=\{x\in dom(g): g(x)\in ran(g)\cap dom(f) \}$ such that $(f\circ g)(x)=f(g(x))$ for all
$x\in dom(g)$, that is
$$ f\circ g: dom(f\circ g)\rightarrow B.$$

Also, we will use the standard definition of identity function: Let
$Id_{A}:A\rightarrow A$ a $\sigma$-function, then we say that
$Id_{A}$ is the \textbf{identity $\sigma$-function} of $A$ iff for
all $x\in A$ we have that $Id_{A}(x)=x$. Now, if we consider the $\sigma$-antiset we have the following definition.

\begin{definition}\label{def 3 antidentity function}
Let $A$ be a $\sigma$-set. If there exists $A^{\st}$ the
$\sigma$-antiset of $A$, then we say that $Id^{\st}_{A}:A\rightarrow
A^{\st}$ is the \textbf{antidentity $\sigma$-function} of $A$ iff for
all $x\in A$ we have that $Id^{\st}_{A}(x)=x^{\st}$.
\end{definition}

We observe that the antidentity $\sigma$-function of $A$ is well
defined. In fact, we consider $x\in A$ then by Theorem \ref{theo 1
antielement is unique} there exists a unique $x^{\st}$ such that
$Id^{\st}_{A}(x)=x^{\st}\in A^{\st}$.

\begin{theorem}\label{theo 3 antifunction of antidentity function}
Let $A$ be a $\sigma$-set such that there exists $A^{\st}$
$\sigma$-antiset of $A$. Then the $\sigma$-antifunction of
$Id^{\st}_{A}$ is the identity $\sigma$-function of $A$, that is
$(Id^{\st}_{A})^{\st}=Id_{A}$.
\end{theorem}

\begin{proof}
Consider $Id^{\st}_{A}:A\rightarrow A^{\st}$ the antidentity
$\sigma$-function of $A$ and $Id_{A}:A\rightarrow A$ the identity
$\sigma$-function of $A$. Now, let $x\in A$, then
$\{Id^{\st}_{A}(x)\}\cup\{Id_{A}(x)\}=\{x^{\st}\}\cup\{x\}=\emptyset$.
Therefore, the $\sigma$-antifunction of $Id^{\st}_{A}$ is the
identity $\sigma$-function of $A$.
\end{proof}

We observe that the $\sigma$-antifunction of $Id_{A}$ is the
antidentity $\sigma$-function of $A$, that is
$Id^{\st}_{A}=(Id_{A})^{\st}$.

Now, suppose that there exists the $\sigma$-antifuncition $f^{\st}$ of a $\sigma$-function $f$ and $f$ has a property $p$, then we will study in that case whether the $\sigma$-antifunction satisfies this property or a similar one.

\begin{lemma}\label{lemma 3.12}
Let $A$ be a $\sigma$-set. If there exists $A^{\st}$ the $\sigma$-antiset of $A$, then for all $B\in 2^{A}$ there exists $D\in 2^{A^{\st}}$ such that $D$ is the $\sigma$-antiset of $B$ (i.e. $D=B^{\st}$).
\end{lemma}

\begin{proof}
Consider $A$ a $\sigma$-set, $A^{\st}$ the $\sigma$-antiset of $A$ and $B\in 2^{A}$. Now, we define the $\sigma$-set $D=\{x\in A^{\st}: x^{\st}\in B\}$. It is clear that $D\in 2^{A^{\st}}$, now we will prove that $D=B^{\st}$. By Theorem \ref{theo 1 antiset is unique} we only prove that $B\cup D=\emptyset$. It is clear that in order to prove that $B\cup D=\emptyset$ it is enough to prove that

\begin{description}
	\item[(a)] $B=B\capst D=\{x\in B: x^{\st}\in D\}.$
	\item[(b)] $D=D\capst B=\{x\in D: x^{\st}\in B\}.$	
\end{description}
Also, we have that $B\capst D\subseteq B$ and $D\capst B\subseteq D$.

(a) Let $x\in B$, then $x\in A$ and in consequence $x^{\st}\in A^{\st}$. Now, since $(x^{\st})^{\st}=x\in B$ then $x^{\st}\in D$. Therefore $x\in B\capst D$ and so $B\capst D= B$.

(b) Let $x\in D$, then $x\in A^{\st}$ and $x^{\st}\in B$. Therefore it is clear that $x\in D\capst B$ and so $D\capst B= D$.

Finally, we obtain that $B\cup D=\emptyset$.
\end{proof}

We observe that by Theorem \ref{theo 3 ex anti ran(f) then  ex antifuntion} and Lemma \ref{lemma 3.12} if $f:A\rightarrow B$ is a $\sigma$-function such that there exists $B^{\st}$ the $\sigma$-antiset of $B$, then there exists $f^{\st}:A\rightarrow B^{\st}$ the $\sigma$-antifunction of $f$.

\begin{theorem}\label{theo 3 func biyect iff antifunc biyect}
Let $A$ and $B$ be $\sigma$-sets such that there exists $B^{\st}$ the $\sigma$-antiset of $B$, $f:A\rightarrow B$ a $\sigma$-function and $f^{\st}:A\rightarrow B^{\st}$ the $\sigma$-antifunction of $f$. Then the following statements hold:
\begin{description}
  \item[(a)] $f$ is a $\sigma$-function onto $B$ if and only if $f^{\st}$ is a $\sigma$-function onto $B^{\st}$.
  \item[(b)] $f$ is a one-one $\sigma$-function if and only if $f^{\st}$ is a $\sigma$-function one-one.
  \item[(c)] $f$ is a bijective $\sigma$-function if and only if $f^{\st}$ is a $\sigma$-function bijective.
\end{description}
\end{theorem}

\begin{proof}
Let $f:A\rightarrow B$ a $\sigma$-function such that there exists $B^{\st}$ the $\sigma$-antiset of $B$ and $f^{\st}:A\rightarrow B^{\st}$ the $\sigma$-antifunction of $f$.

(a) ($\rightarrow$) Suppose that $ran(f)=B$. Since $f^{\st}$ is the $\sigma$-antifunction of $f$, then by Theorem \ref{theo 3 ran(func) antset of ran(antifunc)} $ran(f^{\st})=ran(f)^{\st}=B^{\st}$.

($\leftarrow$) Suppose that $ran(f^{\st})=B^{\st}$. Since $f$ is the $\sigma$-antifunction of $f^{\st}$, then by Theorem \ref{theo 3 ran(func) antset of ran(antifunc)} $ran(f)=ran(f^{\st})^{\st}=(B^{\st})^{\st}=B$.

(b) ($\rightarrow$) Suppose that $f$ is a one-one $\sigma$-function. Let $x,y\in A$ such that $f^{\st}(x)= f^{\st}(y)$ and $x\neq y$. Since $f$ is one-one then $f(x)\neq f(y)$. Now, as $f^{\st}$ is the $\sigma$-antifunction of $f$, then $\{f(x)\}\cup\{f^{\st}(x)\}=\emptyset$ and $\{f(y)\}\cup\{f^{\st}(y)\}=\emptyset$. Therefore, if we define $z=f^{\st}(x)= f^{\st}(y)$ then we obtain two different $\sigma$-antielements of $z$, which is a contradiction by Theorem \ref{theo 1 antielement is unique}.

($\leftarrow$) This proof is similar to the previous one.

(c) This proof is a direct consequence of (a) and (b).
\end{proof}

\begin{figure}
\centering
\includegraphics[width=0.7\textwidth]{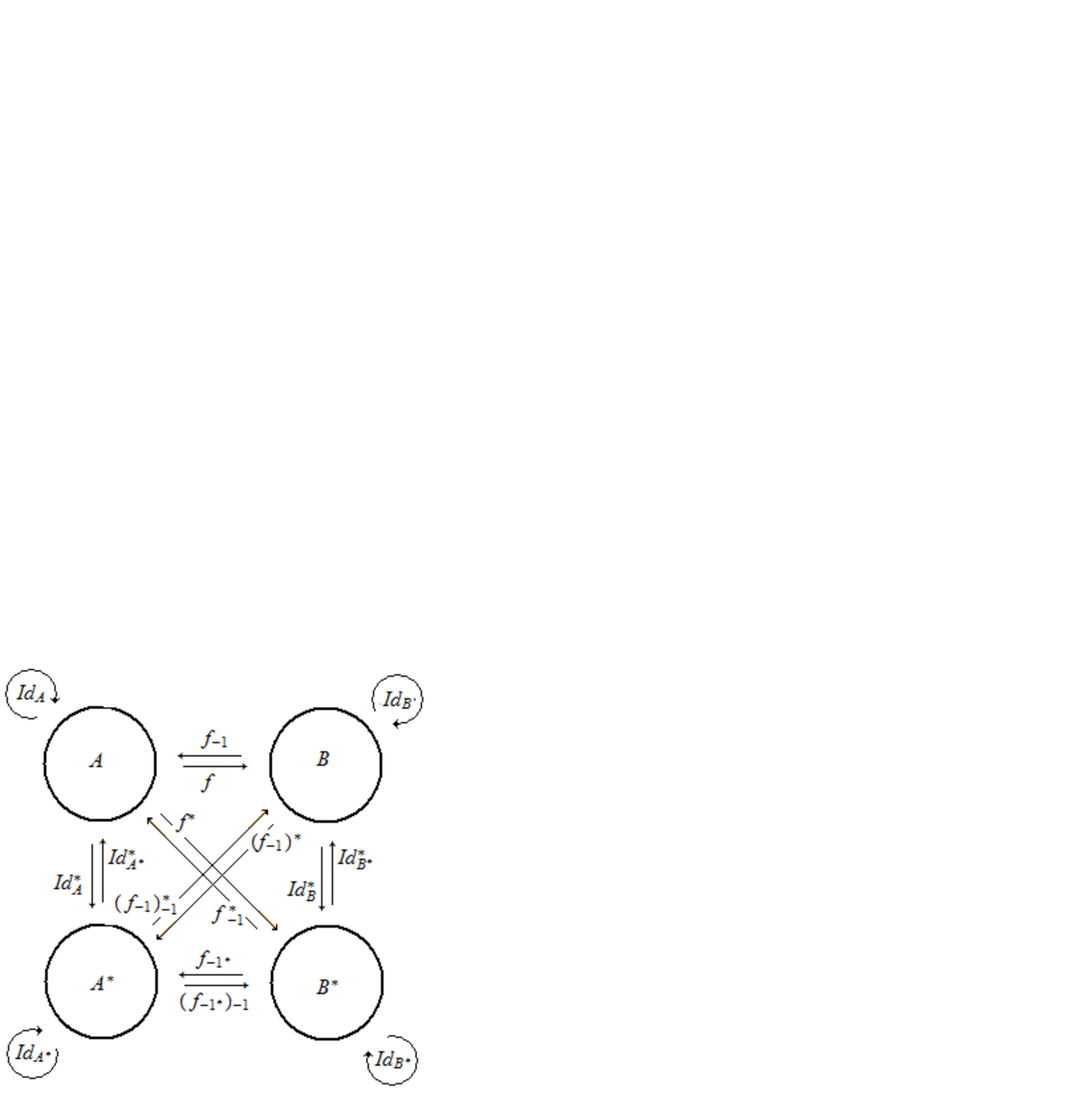}
\caption{Diagram of $\sigma$-functions}
\label{diagram of functions}
\end{figure}

\newpage 

\begin{corollary}
Let $A$ and $B$ be $\sigma$-sets such that there exists $B^{\st}$ the $\sigma$-antiset of $B$, $f:A\rightarrow B$ a $\sigma$-function and $f^{\st}:A\rightarrow B^{\st}$ the $\sigma$-antifunction of $f$. Then there exists $f_{-1}$ the inverse $\sigma$-function of $f$ if and only if there exists $f^{\st}_{-1}$ the inverse $\sigma$-function of $f^{\st}$.
\end{corollary}

\begin{proof}
This proof is a direct consequence of Theorem \ref{theo 3 func biyect iff antifunc biyect}.
\end{proof}

\begin{theorem}\label{theo 3 existence of antinverse}
Let $A$ and $B$ be $\sigma$-sets such that there exist $A^{\st}$ and $B^{\st}$ the $\sigma$-antisets of $A$ and $B$ and $f:A\rightarrow B$ a bijective $\sigma$-function. Then there exists a unique $\sigma$-function $f_{-1^{\st}}:B^{\st}\rightarrow A^{\st}$ (called antinverse of $f$) such that for all $x\in A$ we have that $f_{-1^{\st}}(f^{\st}(x))=x^{\st}$, where $f^{\st}$ is the $\sigma$-antifunction of $f$.	
\end{theorem}

\begin{proof}
\textbf{Existence:} Consider $A$ and $B$, $\sigma$-sets such that there exist $A^{\st}$ and $B^{\st}$ the $\sigma$-antisets of $A$ and $B$ and $f:A\rightarrow B$ a bijective $\sigma$-function. By Theorems \ref{theo 3 ex anti ran(f) then  ex antifuntion} and \ref{theo 3 func biyect iff antifunc biyect} there exists $f^{\st}:A\rightarrow B^{\st}$ the $\sigma$-antifunction of $f$ and it is bijective. Therefore, there exists $f^{\st}_{-1}:B^{\st}\rightarrow A$ the inverse of $f^{\st}$. Since there exists $A^{\st}$ the $\sigma$-antiset of $A$ by Theorem \ref{theo 3 ex anti ran(f) then  ex antifuntion} there exists $(f^{\st}_{-1})^{\st}:B^{\st}\rightarrow A^{\st}$ the $\sigma$-antifunction of $f^{\st}_{-1}$. Therefore, we define $f_{-1^{\st}}=(f^{\st}_{-1})^{\st}$. Since $ran(f^{\st})=dom(f_{-1^{\st}})=B^{\st}$ then we can define $f_{-1^{\st}}\circ f^{\st}: A\rightarrow A^{\st}$. Now, let $x\in A$ then it is clear that $f^{\st}(x)\in B^{\st}$ and $f^{\st}_{-1}(f^{\st}(x))=x$ because $f^{\st}_{-1}$ is the inverse $\sigma$-function of $f^{\st}$. Since $(f^{\st}_{-1})^{\st}$ is the $\sigma$-antifunction of $f^{\st}_{-1}$ then $\{f^{\st}_{-1}(f^{\st}(x))\}\cup\{(f^{\st}_{-1})^{\st}(f^{\st}(x))\}=\emptyset$. Therefore $\{x\}\cup\{f_{-1^{\st}}(f^{\st}(x))\}=\emptyset$ and so $f_{-1^{\st}}(f^{\st}(x))=x^{\st}$.

\textbf{Uniqueness:} Suppose that there exists $\hat{f}:B^{\st}\rightarrow A^{\st}$ a $\sigma$-function such that for all $x\in A$ we have that $\hat{f}(f^{\st}(x))=x^{\st}$ and $f_{-1^{\st}}\neq \hat{f}$. Since $f_{-1^{\st}}\neq \hat{f}$ then there exists $y\in B^{\st}$ such that $f_{-1^{\st}}(y)\neq \hat{f}(y)$. As  $f$ is bijective, by Theorem \ref{theo 3 func biyect iff antifunc biyect} we obtain that $f^{\st}$ is bijective. Therefore, since $y\in B^{\st}$ there exists $x\in A$ such that $f^{\st}(x)=y$. Finally, we obtain two different $\sigma$-antielements of $x$, $f_{-1^{\st}}(y)$ and $\hat{f}(y)$, which is a contradiction.
\end{proof}

We observe that if we consider that $A$ and $B$ are $\sigma$-sets such that there exist $A^{\st}$ and $B^{\st}$ the $\sigma$-antisets of $A$ and $B$ and $f:A\rightarrow B$ a bijective $\sigma$-function, then we can obtain 16 different $\sigma$-functions, that is
\begin{enumerate}
  \item $f:A\rightarrow B$ a bijective $\sigma$-function.
	\item $f_{-1}:B\rightarrow A$ the inverse $\sigma$-function of $f$.
	\item $(f_{-1})^{\st}:B\rightarrow A^{\st}$ the $\sigma$-antifunction of $f_{-1}$.
	\item $(f_{-1})^{\st}_{-1}: A^{\st}\rightarrow B$ the inverse $\sigma$-function of $(f_{-1})^{\st}$.	
	\item $f^{\st}:A\rightarrow B^{\st}$ the $\sigma$-antifunction of $f$.
	\item $f^{\st}_{-1}:B^{\st}\rightarrow A$ the inverse $\sigma$-function of $f^{\st}$.
	\item $f_{-1^{\st}}:B^{\st}\rightarrow A^{\st}$ the antinverse $\sigma$-function of $f$.
	\item $(f_{-1^{\st}})_{-1}:A^{\st}\rightarrow B^{\st}$ the inverse $\sigma$-function of $f_{-1^{\st}}$.
	\item $Id_{A}:A\rightarrow A$ the identity $\sigma$-function of $A$.
	\item $Id_{B}:B\rightarrow B$ the identity $\sigma$-function of $B$.
	\item $Id_{A^{\st}}:A^{\st}\rightarrow A^{\st}$ the identity $\sigma$-function of $A^{\st}$.
	\item $Id_{B^{\st}}:B^{\st}\rightarrow B^{\st}$ the identity $\sigma$-function of $B^{\st}$.
	\item $Id^{\st}_{A}:A\rightarrow A^{\st}$ the antidentity $\sigma$-function of $A$.
	\item $Id^{\st}_{A^{\st}}:A^{\st}\rightarrow A$ the antidentity $\sigma$-function of $A^{\st}$.
	\item $Id^{\st}_{B}:B\rightarrow B^{\st}$ the antidentity $\sigma$-function of $B$.
	\item $Id^{\st}_{B^{\st}}:B^{\st}\rightarrow B$ the antidentity $\sigma$-function of $B^{\st}$.
\end{enumerate}

See Figure \ref{diagram of functions}. The properties of this $\sigma$-function will be studied in future works.

\end{section}

%%%%%%%%%%%%%%%%
% bibliography
%%%%%%%%%%%%%%%

% Set bibliography items using the "thebibliography" environment  and following
% the style used by the AMS journals.
%
% If the bibliography is generated by a bibtex database, use "amsplain" or
% "amsalpha" as bibliography style

\end{document}